\documentclass[12pt, leqno]{article}
\usepackage{amsmath}
\usepackage{amsfonts}
\usepackage{amssymb}
\usepackage{amsthm}
\usepackage{graphics}
\usepackage{graphicx}

\newtheorem{theorem}{Theorem}[section]

\newcommand{\be}{\begin{equation}}
\newcommand{\ee}{\end{equation}}
\begin{document}
\begin{center}
{\bf EXTENDED SIBUYA DISTRIBUTION IN THE SUBCRITICAL MARKOV BRANCHING PROCESSES  }

\vspace{2cm}
Penka Mayster, Assen Tchorbadjieff

\end{center}

\vspace{1.5cm}
\begin{abstract}

The subcritical Markov branching  process $X(t)$ starting with one particle as the initial condition has the ultimate extinction probability $q = 1$. The branching mechanism in consideration is defined by the mixture of logarithmic distributions on the  nonnegative integers.  The purpose of the present paper is to prove that in this case the random number of particles $X(t)$ alive at time $t>0$ follows the shifted extended Sibuya distribution, with parameters depending on the time $t>0$. The conditional limit probability is exactly the logarithmic series distribution supported by the positive integers.
\end{abstract}

\vspace{1cm} {\bf Keywords:} branching process, mixture of logarithmic distributions, extended Sibuya distribution, conditional limit probability


\section{Introduction}

The Sibuya distribution is defined by the probability generating function (p.g.f.)
$$h(s)=1-(1-s)^\gamma,\quad 0<\gamma<1,\quad h(0)=0,\quad h'(1)=\infty. $$
It appears as a particular case of digamma law \cite{SI}.
Then generalised and scaled Sibuya distribution
is considered in \cite{KO}. An overview of the Sibuya-like distributions is given in \cite{KS}.
Extended Sibuya distribution with finite mean \cite{GL} is the natural exponential family extension of the Sibuya distribution with p.g.f.
$$h(s)=\frac{1-(1-b s)^\gamma}{1-(1-b)^\gamma},\quad 0<\gamma<1,\quad 0<b<1,\quad h'(1)<\infty. $$

Branching processes are the fundamental models describing the particles' reproduction development in biology and physics, as introduced in \cite{AN, H, SE}.
The subcritical Markov branching process (MBP) yields the remarkable property, having certain ultimate extinction, but positive conditional limit law over the non-extinction to the time $t>0$.

The explicit solution to the Kolmogorov equation is well-known for the linear birth-death process. The MBP $X(t), t>0, X(0)=1, $ with geometric reproduction of particles is developed in \cite{TMT,TM}. The p.g.f. $ F(t,s):=E[s^{X(t)}]$ is defined there by the composition of special functions such as the Wright and Lambert-W functions.

In the present short communication, we consider only the subcritical case. The reproduction of particles is defined by the mixture of logarithmic distributions and Dirac measure at zero. In this special model, we find the explicit representation of the distribution for the number of particles alive $X(t), t>0,$ as the shifted extended Sibuya distribution.  The conditional limit probability is exactly the logarithmic series, \cite{FCW, JKK}, distribution supported by positive integers $ N=\{1, 2,  ...\}$. The correspondence between the branching mechanism and conditional limit law is noted in the conclusion of the present communication.

The L\'{e}vy process generated by the infinitely divisible logarithmic distribution, \cite{JKK, SH}, supported by nonnegative integers $ Z_{+}=\{0, 1, 2, . . .\}$ is studied in \cite{MTL}.
\section{Kolmogorov equations}
The MBP $X(t), t>0, X(0)=1,$ is defined by the random lifetime of particles exponentially distributed with constant parameter $K>0$ and reproduction law of the newborn particles.
The main assumption is the independence of the evolution of particles. The offspring number is a random variable $\eta$ with the following probability
mass function (p.m.f.),
$$ \textbf{P}(\eta=0)=\alpha,\quad \textbf{P}(\eta=1)=1-\alpha^2\left(1+\frac{1}{A}\right),\quad A=\sum^\infty_{n=1}\frac{\alpha^n}{n}>0 ,$$
and  for $n=2, 3,...   , $  as follows
$$ \textbf{P}(\eta=n)= \frac{\alpha}{A}\frac{\alpha^n}{n(n-1)}=\frac{\alpha}{A}\left\{\frac{\alpha^n}{n-1}-\frac{\alpha^n}{n}\right\}. $$
Respectively, the probability generating function of reproduction is
\begin{equation}
h(s)=s+\alpha(1-\alpha s)\left\{ 1+\frac{\log(1-\alpha s)}{A}\right\},\quad |s| \leq1, \label{BM}
\end{equation}
under the following restriction on parameters, insuring $0< P(\eta=1)<1$,
\begin{equation}
A=-\log(1-\alpha),\quad 0<\alpha<\alpha_\ast\simeq 0,772638,\quad 1+\frac{1}{A}=\frac{1}{(\alpha_\ast)^2}.\label{rest}
\end{equation}
The first derivative of the reproduction p.g.f. is calculated as
 $$h'(s)=1- \alpha^2-\frac{\alpha^2}{A}+\frac{\alpha^2}{A}\sum^\infty_{n=2}\frac{(\alpha s)^{n-1}}{(n-1)},\quad h'(1):=m=1-\frac{\alpha^2}{A}<1,\quad m>0. $$
 The parameter $m:=E[\eta]$ defines the mean of the offspring numbers. Consequently, \cite{AN, H, SE}, the reproduction is subcritical, and the ultimate extinction probability
 $$q:=\lim_{t\rightarrow \infty} \textbf{P}(X(t)=0)=1.$$
The infinitesimal generating function is given by
 $$ f(s)=K(h(s)-s)= \frac{K \alpha}{A}(1-\alpha s)(A+\log(1-\alpha s)), \quad f'(1)=-\frac{K\alpha^2}{A}.$$
The p.g.f. of the number of particles alive at the positive time $t>0$
\begin{equation}
F(t,s)=\sum^\infty_{k=0}s^k P(X(t)=k|X(0)=1),\quad |s|\leq 1,\label{F}
\end{equation}
yields the backward Kolmogorov equation
\begin{equation}
 \frac{\partial}{\partial t}\left(F(t,s)\right)=f\left(F(t,s)\right),\quad F(0,s)=s. \label{CKs}
 \end{equation}
 The mathematical expectation of the number of particles alive at the positive time $ t>0 $ in the subcritical case has the following exponential decreasing behaviour,
 \begin{equation}
  \textbf{E}[X(t)]:=M(t)=\exp\{f'(1)t\}<1, \quad M'(t)=f'(1)M(t),\quad f'(1) <0 .\label{ME}
 \end{equation}
Equation (\ref{CKs}) is nonlinear, of the type separate differentials, due to the time-homogeneous property.
It is equivalent to
\begin{equation}
\frac{f'(1) dx}{f(x)}=f'(1)dt,\quad x=F(t,s),\quad F(0,s)=s,\quad f'(1)=-\frac{K\alpha^2}{A} . \label{Ks}
 \end{equation}
The infinitesimal generating function is very convenient to calculate the primitive using the following representation
$$\frac{f'(1) dx}{f(x)}=\left(\frac{K\alpha}{A}\right)\left(\frac{-\alpha dx}{f(x)}\right)=\frac{d (A+\log(1-\alpha x))}{A+\log(1-\alpha x)}. $$
We calculate the integral of (\ref{Ks}) in the explicit form
 $$\int^s_{x=0}\frac{f'(1) dx}{ f(x)}=\int^s_{x=0}\frac{ d (A+\log(1-\alpha x))} {A+\log(1-\alpha x)}=\log(A+\log(1-\alpha s))-\log(A).$$
  The implicit solution of equations (\ref{CKs}) and (\ref{Ks}) with the initial condition, $F(0,s)=s $, is written
  $$\log\left(1+\frac{\log(1-\alpha F(t,s))}{A}\right)=\log(e^{f'(1)t})+\log\left(1+\frac{\log(1-\alpha s)}{A}\right). $$
 It can be expressed by the mathematical expectation (\ref{ME}), $ M(t)=e^{f'(1)t}$, as
 \begin{equation}
\log\left(1+\frac{\log(1-\alpha F(t,s))}{A}\right)=\log \left\{M(t) \left(1+\frac{\log(1-\alpha s)}{A}\right)\right\}.\label{Impl}
\end{equation}
\section{The main result}
 \begin{theorem}
Let $X(t), t>0,$ be a subcritical time-homogeneous MBP with a branching mechanism given by the p.g.f. (\ref{BM}) under the restriction (\ref{rest}).
The explicit solution for (\ref{F}) to the backward Kolmogorov equation (\ref{CKs}) is written by the mathematical expectation (\ref{ME})
 with $ M(0)=1, F(0,s)=s,$ as follows
\begin{equation}
F(t,s)= \frac{1}{\alpha}\left\{1-(1-\alpha)\left(\frac{1-\alpha s}{1-\alpha}\right)^{M(t)} \right\},\quad F(t,1)=1,\quad |s|\leq 1 .\label{ES}
\end{equation}
 \end{theorem}
  \begin{proof}
Starting with the implicit solution (\ref {Impl}), we remark only that the logarithmic function is strictly increasing and we transform (\ref {Impl}) into
$$ A+\log(1-\alpha F(t,s))=  M(t)(A+\log(1-\alpha s)).$$
  Remembering that $ A=-\log(1-\alpha)$ we write
  $$ \log\left(\frac{1-\alpha F(t,s)}{1-\alpha}\right) =  \log\left\{\left(\frac{1-\alpha s}{1-\alpha}\right)^{M(t)}\right\}.$$
Obviously, the representation (\ref{ES}) is equivalent to the following
 $$ \frac{1-\alpha F(t,s)}{1-\alpha} = \left(\frac{1-\alpha s}{1-\alpha}\right)^{M(t)},\quad M(t)=e^{f'(1)t},\quad M(0)=1. $$
    \end{proof}
  The first derivative of $ F(t,s)$ (\ref{ES}) in $t>0$ is a derivative of the exponential function
 $$F'_{t}(t,s)=-\left(\frac{1-\alpha}{\alpha}\right) \left( \frac{1-\alpha s}{1-\alpha}\right)^{M(t)}\log \left( \frac{1-\alpha s}{1-\alpha}\right)M'(t),\quad M'(t)=f'(1) M(t),$$
and we have
$$ F'_{t}(t,s)=\frac{K\alpha (1-\alpha)}{A}\left( \frac{1-\alpha s}{1-\alpha}\right)^{M(t)} \log \left\{ \left( \frac{1-\alpha s}{1-\alpha}\right)^{M(t)}\right\}. $$
The first derivative of $ F(t,s)$ (\ref{ES}) in $s>0$ is a derivative of the power function
$$F'_{s}(t,s)=M(t) \left( \frac{1-\alpha s}{1-\alpha}\right)^{M(t)-1}. $$
Obviously, $F'_{t}(t,s) $ and $F'_{t}(t,s) $ fulfill the backward and forward Kolmogorov equations,
$$F'_t(t,s)=f(F(t,s)),\quad F'_t(t,s)=f(s) F'_s (t,s) ,\quad F(0,s)=s,\quad F(t,1)=1.$$
The extinction probability to the positive time $t>0$ is equal to
 $$\textbf{P}(X(t)=0)=F(t,0)=\frac{1}{\alpha}\left\{1-\frac{1-\alpha}{(1-\alpha)^{M(t)}}\right\}. $$
Respectively, the survival probability to $ t>0$ is written  as follows,
\begin{equation}\label{sur}
1-F(t,0)=\left(\frac{1-\alpha}{\alpha}\right)\left\{\frac{1}{(1-\alpha)^{M(t)}}-1\right\}. \end{equation}
The ultimate extinction probability is given by the following limits
$$q:=\lim_{t\rightarrow\infty}F(t,0)=\frac{1-(1-\alpha)}{\alpha}=1,\quad \lim_{t\rightarrow\infty}M(t)=\lim_{t\rightarrow\infty}\exp\left(\frac{-K\alpha^2 t}{A}\right)=0.$$
In order to define the consecutive derivatives, we introduce the notation of the falling factorials as follows,
$$[x]_{n \downarrow}=x(x-1)...(x-n+1)= \frac{\Gamma(x+1)}{\Gamma (x+1-n)}.$$
The consecutive derivatives of the p.g.f. $F(t,s)$  in $s$ are written as follows
 \begin{equation} \label{pgf}
 F^{(n)}_s(t,s)=-\left(\frac{1-\alpha}{\alpha}\right) \frac{(1-\alpha s)^{M(t)-n}\alpha^n (-1)^n[M(t)]_{n\downarrow}}{(1-\alpha)^{M(t)}}.
\end{equation}
The values $F^{(n)}_s(t,0) $ and $F^{(n)}_s(t,1) $ of (\ref{pgf}) express the p.m.f. and factorial moments in the following theorem.
\begin{theorem}
Let $X(t), t>0,$ be a subcritical time-homogeneous MBP with a branching mechanism given by the p.g.f. (\ref{BM}) under the restriction (\ref{rest}). Then the p.m.f. and the factorial moments of the number of particles alive to the time $t>0$
are given by the mathematical expectation $M(t)$ (\ref{ME}) as follows
 \begin{equation} \label{pmf}
\textbf{P}(X(t)=n):=\frac{F^{(n)}_s(t,0)}{n!}=\frac{-1}{n!}\left(\frac{1-\alpha}{\alpha}\right) \frac{\alpha^n (-1)^n[M(t)]_{n\downarrow}}{(1-\alpha)^{M(t)}}>0,
 \end{equation}
 and
 \begin{equation} \label{ffm}
\textbf{E}[X(t)]_{n\downarrow}:=F^{(n)}_s(t,1))=-\left(\frac{1-\alpha}{\alpha}\right) \frac{\alpha^n (-1)^n[M(t)]_{n\downarrow}}{(1-\alpha)^{n}}>0.
 \end{equation}
 \end{theorem}

 \subsection{The conditional p.m.f. and conditional factorial moments}
 The conditional p.m.f. for $ n=1,2,...$, knowing the non-extinction (\ref{sur}) and p.m.f. (\ref{pmf}), is
$$
  \textbf{P}(X(t)=n|X(t)>0)=\frac{-1}{n!}\left(\frac{1-\alpha}{\alpha}\right) \frac{\alpha^n (-1)^n[M(t)]_{n\downarrow}}{(1-\alpha)^{M(t)}}\frac{1}{\left(\frac{1-\alpha}{\alpha}\right)\left\{\frac{1}{(1-\alpha)^{M(t)}}-1\right\}},
$$
 and equivalently
  \begin{equation} \label{cmf}
  \textbf{P}(X(t)=n|X(t)>0)=\frac{-1}{n!}\left( \frac{\alpha^n (-1)^n[M(t)]_{n\downarrow}}{1-(1-\alpha)^{M(t)}}\right).
 \end{equation}
 The conditional factorial moments for $ n=1,2,...$, knowing the non-extinction (\ref{sur}) and factorial moments (\ref{ffm}), is
$$  \textbf{E}([X(t)]_{n\downarrow}|X(t)>0)=-\left(\frac{1-\alpha}{\alpha}\right) \frac{\alpha^n (-1)^n[M(t)]_{n\downarrow}}{(1-\alpha)^{n}}\frac{1}{\left(\frac{1-\alpha}{\alpha}\right)\left\{\frac{1}{(1-\alpha)^{M(t)}}-1\right\}},
$$
 and equivalently
  \begin{equation} \label{cfm}
  \textbf{E}([X(t)]_{n\downarrow}|X(t)>0)=-\left( \frac{\alpha^n (-1)^n[M(t)]_{n\downarrow}}{1-(1-\alpha)^{M(t)}}\right)\left(\frac{(1-\alpha)^{M(t)}}{(1-\alpha)^{n}}\right).
 \end{equation}
\subsection{Limit theorem}
The p.g.f. of the conditional probability is written as follows
$$\sum^\infty_{n=1}s^n \textbf{P}(X(t)=n|X(t)>0)=\frac{F(t,s)-F(t,0)}{1-F(t,0)}=\frac{1-(1-\alpha s)^{M(t)}}{1-(1-\alpha)^{M(t)}}. $$
 The asymptotic behaviour of the conditional p.m.f. and factorial moments when $ t\rightarrow \infty$ is based on the following equivalence for the exponential function in the neighborhood of zero, as $ \lim_{t\rightarrow \infty}M(t)=0$,
 $$(1-\alpha s)^{M(t)}=e^{(\log(1-\alpha s))M(t)}\sim 1+(\log(1-\alpha s))M(t).$$
 Knowing the previous equivalence
we obtain the p.g.f. of the positive random variable $ \xi \subset (1,2,...),$ (under the conditional limit) as follows
 $$F^\ast(s):=\textbf{E}[s^\xi]= \lim_{t\rightarrow \infty}\frac{1-(1-\alpha s)^{M(t)}}{1-(1-\alpha)^{M(t)}},$$
in its explicit form
$$F_\ast(s)=\frac{-\log(1-\alpha s)}{-\log(1-\alpha )}=\frac{1}{A}\sum^\infty_{k=1}\frac{(\alpha s)^k}{k},\quad A=-\log(1-\alpha)>0 .$$
  Having in mind the definition of decreasing factorials and decomposition
  $$ [M(t)]_{n\downarrow}=M(t)[M(t)-1]_{(n-1)\downarrow},\quad M(t)<1, $$
we write the limit of (\ref{cmf}) and (\ref{cfm}) as follows
 $$\lim_{t\rightarrow \infty}\textbf{P}(X(t)=n|X(t)>0)=\frac{\alpha^n(n-1)!}{A n!} =\frac{\alpha^n}{A n}, \quad n=1,2,...,$$
and
 $$\lim_{t\rightarrow \infty} \textbf{E}([X(t)]_{n\downarrow}|X(t)>0)=\frac{\alpha^n(n-1)!}{A (1-\alpha)^n} = \frac{(n-1)!}{A}\left(\frac{\alpha}{1-\alpha}\right)^n, \quad n=1,2,....$$
 Each one of these limits implies the two others.

Finally, the important property  for the behavior of the process, like the conditional limit
probability is summarised in the following theorem.
\begin{theorem} \label{tm:gen.f}
Let $X(t), t>0,$ be a subcritical time-homogeneous MBP with a branching mechanism given by the p.g.f. (\ref{BM}) and restriction (\ref{rest}). Then the conditional limit probability exists and is given by the random variable $\xi $,
$$ \lim_{t\rightarrow \infty} \textbf{P}(X(t)=n|X(t>0))=p_n=P(\xi=n) , \quad n=1,2,...,$$
with p.g.f. $$ F_\ast(s)= \sum^\infty_{n=1} p_n s^n=\frac{-\log(1-\alpha s)}{A}, \quad A=\sum^\infty_{n=1}\frac{\alpha^n}{n}, \quad |s|\leq 1, $$
and
$$\textbf{P}(\xi=n)=\frac{\alpha^n}{An},\quad \textbf{E}[\xi]_{n\downarrow}=\frac{(n-1)!}{A } \left(\frac{\alpha}{1-\alpha}\right)^{n}>0.$$
\end{theorem}
\section{Conclusion remark}
The explicit form of the p.g.f. $F(t, s)$ (\ref{ES}) can be classified as Sibuya-like shifted extended distribution \cite{KS}.
The correspondence between the branching mechanism and asymptotic
behavior is a very interesting problem for future consideration. In the
following Table, we give only several examples, see \cite{TMT,TM} for the geometric reproduction of particles.
\begin{table}[t]
	\caption{examples for the correspondence - p.g.f. $h(s)$ and $F_\ast(s)$}\label{t1}
	\begin{tabular}{@{}ll}
		\hline
		$ h(s),\quad h'(1)=m<1,$ &$ F_\ast(s)$ \\
		&  \\
		$ h(s)=s+\alpha(1-\alpha s)\left\{ 1+\frac{\log(1-\alpha s)}{A}\right\}$  & $F_\ast(s)= \frac{-\log(1-\alpha s)}{A},\quad A=-\log(1-\alpha) $\\
		&  \\
		$h(s)=\frac{1}{1+m-ms}$   & $  F_\ast(s)=1-\frac{1-s}{(1-m s)^m} $\\
&  \\
$h(s)=1+\frac{m}{2}(s^2-1) ,\quad m=\frac{2\varrho}{1+\varrho}$   & $ F_\ast(s)=\frac{(1-\varrho)s}{1-\varrho s},\quad  \varrho=\frac{m}{2-m}$\\
		&  \\
		$h(s)=1-m +ms$& $F_\ast(s)=s $ \\	
		&  \\
		\hline
	\end{tabular}
\end{table}

\noindent
Assen Tchorbadjieff, Penka Mayster\\
Institute of Mathematics and Informatics, \\
Bulgarian Academy of Sciences,\\
Acad. G. Bonchev street, Bloc 8, 1113 Sofia. \\
Email address: atchorbadjieff@math.bas.bg\\
Email address: penka.mayster@math.bas.bg

\end{document}